\documentclass{article}
\usepackage{latexsym,amsfonts,amsmath,amsthm,makeidx}
\usepackage{makeidx}

\makeindex
\newtheorem{definition}{Definition}[section]
\newtheorem{theorem}{Theorem}[section]
\newtheorem{lemma}{Lemma}[section]
\newtheorem{corollary}{Corollary}[section]
\newtheorem{proposition}{Proposition}[section]
\newtheorem{remark}{Remark}[section]
\newtheorem{example}{Example}[section]
\newcommand{\RN}{\mathbb R^N}
\newcommand{\om}{\Omega}

\newcommand{\iy}{\infty}

\newcommand{\s}{\section}
\newcommand{\dd}{\delta}
\newcommand{\DD}{\Delta}
\newcommand{\g}{\gamma}
\newcommand{\G}{\Gamma}
\newcommand{\na}{\nabla}

\newcommand{\pa}{\partial}
\newcommand{\si}{\sigma}

\newcommand{\R}{\mathbb R}
\newcommand{\al}{\alpha}

\newcommand{\ti}{\tilde}
\newcommand{\bb}{\beta}

\newcommand{\rg}{\rightarrow}

\newcommand{\e}{\varepsilon}
\newcommand{\vp}{\varphi}

\newcommand{\lab}{\label}
\newcommand{\bt}{\begin{theorem}}
\newcommand{\et}{\end{theorem}}
\newcommand{\bl}{\begin{lemma}}
\newcommand{\el}{\end{lemma}}
\newcommand{\bd}{\begin{definition}}
\newcommand{\ed}{\end{definition}}
\newcommand{\bc}{\begin{corollary}}
\newcommand{\ec}{\end{corollary}}
\newcommand{\bp}{\begin{proof}}
\newcommand{\ep}{\end{proof}}
\newcommand{\bx}{\begin{example}}
\newcommand{\ex}{\end{example}}
\newcommand{\bi}{\begin{exercise}}
\newcommand{\ei}{\end{exercise}}
\newcommand{\bo}{\begin{proposition}}
\newcommand{\eo}{\end{proposition}}
\newcommand{\br}{\begin{remark}}
\newcommand{\er}{\end{remark}}
\newcommand{\be}{\begin{equation}}
\newcommand{\ee}{\end{equation}}
\newcommand{\ba}{\begin{align}}
\newcommand{\ea}{\end{align}}
\newcommand{\bn}{\begin{enumerate}}
\newcommand{\en}{\end{enumerate}}
\newcommand{\bg}{\begin{align*}}
\newcommand{\bcs}{\begin{cases}}
\newcommand{\ecs}{\end{cases}}

\newcommand{\bean}{\begin{eqnarray*}}
\newcommand{\eean}{\end{eqnarray*}}

\numberwithin{equation}{section}

\begin{document}

\title{\bf{Standing Waves for nonlinear Schr\"{o}dinger Equations involving critical growth
 \thanks{This work is supported by NSFC(10871109, 11025106, 10771212)}}}
\date{}
\author{{\bf Jianjun Zhang \& Zhijie Chen \& Wenming Zou}\\
\footnotesize {\it  Department of Mathematical Sciences, Tsinghua
University,}\\
\footnotesize {\it Beijing 100084, China}\\
}

\numberwithin{equation}{section}

\maketitle

\vskip0.36in

\begin{center}
\begin{minipage}{120mm}
\begin{center}{\bf Abstract}\end{center}

We consider the following singularly perturbed nonlinear elliptic problem:
$$-\e^2\Delta u+V(x)u=f(u),\ u\in H^1(\mathbb{R^N}),$$
where $N\ge 3$ and the nonlinearity $f$ is of {\it critical growth}. In this paper, we construct a solution $u_\e$ of the above problem which concentrates at an isolated component of positive local minimum points of $V$ as $\e\rg 0$  under certain conditions on $f$.
 Our result completes the study made in  some  very recent works  in the
sense that, in those  papers  only the subcritical growth was considered.\\

\end{minipage}
\end{center}

\vskip0.16in

\s{Introduction}
\renewcommand{\theequation}{1.\arabic{equation}}

In this paper, we shall be concerned with the existence and concentration of positive solutions for the following singular perturbed elliptic problem with critical growth:
\be
\label{q1}
-\e^2\DD v+V(x)v=f(v),\ v>0,\ v\in H^1(\mathbb{R^N}),
\ee
where $N\ge 3$. For $\e>0$ sufficiently small, these standing waves are referred to as semi-classical states. In the sequel, we assume that the potential function $V$ satisfies the following conditions:
\begin{itemize}

\item [(V1)] $V\in C(\RN,\R)$ and $0<V_0:=\inf_{x\in\RN}V(x)$;

\item [(V2)] There is a bounded domain $O$ such that $$m:\equiv\inf_{x\in O}V(x)<\min_{x\in \partial O}V(x).$$

\end{itemize}

In 2007, Byeon and Jeanjean \cite{byeon} considered the concentration phenomenon of the above  problem (\ref{q1}) and developed a new variational method to explore what are the essential features which guarantee the existence of localized ground states. The considered the following conditions:
\begin{itemize}

\item [$({f_1})$] $f\in C(\R, \R)$ such that $f(t)=0$ for $t\leq 0$ and $\lim_{t\to 0}f(t)/t=0;$

\item [$({f_2})$]  there exists $p\in (1, (N+2)/(N-2))$ such that $\lim\sup_{t\to \infty}f(t)/t^p<\infty;$

\item [$({f_3})$]   there exists $T>0$ such that $\frac{m}{2}T^2<F(T) :\equiv\int_0^Tf(t)dt.$

\end{itemize}

Let $$\mathcal{M}\equiv\{x\in O: V(x)=m\}.$$

\noindent{\bf Theorem A} (see \cite{byeon}) {\it Suppose that $(V1)$-$(V2)$ and $(f_1)$-$(f_3)$.Then for sufficiently small $\e>0$, $(\ref{q1})$ admits a positive solution $v_{\e}$, which satisfies
\begin{itemize}
\item [(i)] there exists a maximum point $x_\e$of $v_\e$ such that $\lim_{\e\rg 0}dist(x_\e,\mathcal{M})=0$ and for any such $x_\e$,$w_\e(x)\equiv v_\e(\e x+x_\e)$ converges (up to a subsequence) uniformly to a least energy solution of
\be\lab{BB}
-\DD u+mu=f(u),\ \ u>0,\ \ u\in H^1(\RN),
\ee
\item [(ii)] $v_\e(x)\le C\exp(-\frac{c}{\e}|x-x_\e|)$ for some $c,C>0$.
\end{itemize}}
In \cite{byeon}, Byeon and Jeanjean believed that $(f_1)$-$(f_3)$ are almost optimal for the subcritical case. Hypotheses $(f_1)$-$(f_3)$ are called Berestycki-Lions\ conditions, which were firstly  proposed in a classical paper \cite{Lions} to guarantee the existence of ground states of (\ref{BB}) in the subcritical case. It follows from Pohozaev's identity (cf.\cite{Pohozaev}), that $(f_3)$ is necessary and that for $f(u)=u^p$ with $p\ge \frac{N+2}{N-2}$, there exists no nontrivial solutions  in $H^1(\RN)$. Thus,  Berestycki-Lions\ conditions  are almost optimal for the existence of solutions for (\ref{BB}) (cf. \cite{byeon}).

\vskip0.1in

%%%%%%%%%%%%%%%%%%%%%%%%%%%%%%%%%%%%%%%%%%%%%%%%%%%%%%%%

Since in $(f_2)$ above,   $p\in (1, (N+2)/(N-2))$ characteristics the problem to be of subcritical growth. A natural open problem which has not been settled before the case of critical growth, is whether the results like Theorem  A  hold if $f$ is of critical growth?  The purpose of this paper is to complete the study for such an open problem  with critical exponent growth.
Before making more comments on the background of such singularly  perturbed nonlinear elliptic problems, we state  the main result of this paper first.  It is well known, the critical exponent growth makes the problem very tough, more assumptions are of course needed. We now assume that $f\in C(\R,\R)$ and satisfies:
\begin{itemize}
\item [(F1)] $\lim_{t\rg 0}\frac{f(t)}{t}=0$.
\item [(F2)] $\lim_{t\rg \infty}\frac{f(t)}{t^{\frac{N+2}{N-2}}}=\kappa>0$.
\item [(F3)] There exist $C>0$ and $p<2^{\ast}$ such that $f(t)\ge\kappa t^{\frac{N+2}{N-2}}+Ct^{p-1}$ for $t\ge 0$.
\end{itemize}
The main theorem of this paper reads as
\bt\lab{Theorem 1} Let $p>2,N\ge 4$ or $p>4, N=3$ and suppose that $(V1)$-$(V2)$ and $(F1)$-$(F3)$.
Then for sufficiently small $\e>0$, $(\ref{q1})$ admits a positive solution $v_{\e}$, which satisfies
\begin{itemize}
\item [(i)] there exists a maximum point $x_\e$of $v_\e$ such that $\lim_{\e\rg 0}dist(x_\e,\mathcal{M})=0$ and for any such $x_\e$,$w_\e(x)\equiv v_\e(\e x+x_\e)$ converges (up to a subsequence) uniformly to a least energy solution of
\be\lab{q3}
-\DD u+mu=f(u),\ \ u>0,\ \ u\in H^1(\RN),
\ee
\item [(ii)] $v_\e(x)\le C\exp(-\frac{c}{\e}|x-x_\e|)$ for some $c,C>0$.
\end{itemize}
\et

\br
Without loss of generality, in the present paper we can assume that $V_0=\kappa=1$.
\er

\br To ensure the existence of ground states to $(\ref{q3})$, the assumption $(F3)$ plays a crucial role.
Without $(F3)$, the assumptions $(F1)$-$(F2)$ can not guarantee the existence of ground states of $(\ref{q3})$. We can give a counterexample, i.e., $f(s)=\kappa|s|^{2^{\ast}-2}s$. Then $f$ satisfies the assumptions $(F1)$-$(F2)$ except $(F3)$. But it is easy to verify with the help of Pohoz\v{a}ev's identity that $(\ref{q3})$ has no nontrivial solutions.
\er

In the study of singularly  perturbed problems,   the limit problem (\ref{q3}) plays a crucial role. In   \cite{Jean}, Jeanjean and Tanaka showed  under the Berestycki-Lions  conditions  $(f_1)$-$(f_3)$ that  the subcritical problem (\ref{q3}) exists a least energy solution, which is also a mountain pass solution. Due to the  lack of compact embedding of $H^1(\RN)\hookrightarrow L^{2^{\ast}}(\RN)$, for critical nonlinearity $f$, the existence of ground states of problem (\ref{q3}) becomes rather complicated. Very recently, Alves, Souto and Montenegro \cite{souto} studied  the existence of ground state solutions for problem (\ref{q3}) with critical growth in $\RN (N\ge 2)$. For $N\ge 3$, they assume that $f\in C(\R, \R)$ and satisfies
\begin{itemize}
\item [(G1)] $\lim_{t\rg 0^+}\frac{f(t)}{t}=0$;
\item [(G2)] $\limsup_{t\rg \infty}\frac{f(t)}{t^{2^{\ast}-1}}\le 1$;
\item [(G3)] $2F(t)\le tf(t)$ for all $t\ge 0$;
\item [(G4)] There exist $\lambda>0$ and $2<p<2^{\ast}$ such that $f(t)\ge \lambda t^{p-1}$ for $t\ge 0$.
\end{itemize}
They established the existence of the ground state to (\ref{q3}). But the proof in \cite{souto} strongly depends on large $\lambda$, that is, problem (\ref{q3}) has a ground state if $\lambda>\lambda_0$, where $\lambda_0$ is a positive constant and a complicated explicit formula of $\lambda_0$ is given there. For small $\lambda>0$, it still remains unknown that whether problem (\ref{q3}) has a ground state.  In \cite{Zhang-Zou}, we proved that problem (\ref{q3}) has a ground state with the assumptions $(F_1)$-$(F_3)$. Meanwhile, we show that a ground state of (\ref{q3}) is a mountain pass solution. More properties are also claimed.

\vskip0.1in
%%%%%%%%%%%%%%%%%%%%%%%%%%%%%%%%%%%%%%%%%%%%%%%%%%%%%%%%%%%%%%%%%%%%%%%%%%%%%%%%%%%%%%%%%%%%%%%%%%%%%%%%%%%%%%%%%%%%%%%%%%%%%%%%%%%%%%%%%%%%%%%%%%%
%%%%%%%%%%%%%%%%%%%%%%%%%%%%%%%%%%%%%%%%%%%%%%%%%%%%%%%%%%%%%%%%%%%%%%%%%%%%%%%%%%%%%%%%%%%%%%%%%%%%%%%%%%%%%%%%%%%%%%%%%%%%%%%%%%%%%%%%%%%%%%%%%%%
%%%%%%%%%%%%%%%%%%%%%%%%%%%%%%%%%%%%%%%%%%%%%%%%%%%%%%%%%%%%%%%%%%%%%%%%%%%%%%%%%%%%%%%%%%%%%%%%%%%%%%%%%%%%%%%%%%%%%%%%%%%%%%%%%%%%%%%%%%%%%%%%%%%
%%%%%%%%%%%%%%%%%%%%%%%%%%%%%%%%%%%%%%%%%%%%%%%%%%%%%%%%%%%%%%%%%%%%%%%%%%%%%%%%%%%%%%%%%%%%%%%%%%%%%%%%%%%%%%%%%%%%%%%%%%%%%%%%%%%%%%%%%%%%%%%%%%%
%%%%%%%%%%%%%%%%%%%%%%%%%%%%%%%%%%%%%%%%%%%%%%%%%%%%%%%%%%%%%%%%%%%%%%%%%%%%%%%%%%%%%%%%%%%%%%%%%%%%%%%%%%%%%%%%%%%%%%%%%%%%%%%%%%%%%%%%%%%%%%%%%%%
Now let us say more on the background for problems like $(\ref{q1})$. In recent years, singularly perturbed problems have  been widely studied by many researchers, and related results can been seen in \cite{byeon,byeon2,byeon3,byeon1,byeon4,Felmer,Pino,Ni-Wei,Ni-Takagi1,Ni-Takagi2,ramos}.
By denoting $u(x)=v(\e x)$ and $V_\e(x)=V(\e x)$, $(\ref{q1})$ is equivalent to
\be
\label{q2}
-\DD u+V_\e(x)u=f(u),\ u>0,\ u\in H^1(\mathbb{R^N}).
\ee
An interesting class of solutions of (\ref{q1}) are families of solutions which concentrate and develop spike layers around certain point in $\om$ as $\e\rg 0$. To study the concentration phenomena of solutions for problem (\ref{q1}), the  problem (\ref{q3}) plays an important role which is called the limit problem of (\ref{q2}).

\vskip0.1in

Recall that Floer and Weinstein \cite{F-W} first studied the existence of single peak solutions for $N=1$ and $f(s)=s^3$. They construct a single peak solution which concentrates around any given non-degenerate critical point of $V$. In higher dimension, for $f(s)=|s|^{p-2}s, p\in (2,2^{\ast})$, Oh \cite{Oh} established a similar result as in \cite{F-W}. In \cite{F-W,Oh}, their arguments are based on a Lyapunov-Schmidt reduction, for which they needed to characterize the kernel of the linearized operator $L:=-\DD+V(x_0)-f'(U)$, where $U$ is the ground state of the following autonomous problem: for fixed $x_0\in\RN$,
\be\lab{ap1}
-\DD u+V(x_0)u=f(u),\ \ \mbox{in}\ \RN,\ \ \ v\in H^1(\RN).
\ee
Moreover, they also required some monotonicity condition of nonlinearity $f$ and uniqueness condition of ground states of (\ref{ap1}). Precisely, they assumed that $f\in C^{0,1}(\R,\R)$ and
\begin{itemize}
\item [(H1)] $f(t)/t$ is non-decreasing on $(0,\infty);$
\item [(H2)] there exists a unique radially symmetric solution $U
\in H^{1,2}(\RN)$ for $\Delta u - u + f(u) = 0, u > 0$ in $\RN$
such that if $\Delta V - V + f^\prime(U)V = 0$ and $ V \in
H^{1,2}(\RN),$ then $V = \sum_{i=1}^n a_i\frac{\pa U}{\pa x_i}$
for some $a_1,\cdots,a_n \in \R.$
\end{itemize}
Subsequently, when ground states of (\ref{ap1}) are unique and non-degenerate, Ambrosetti, Badiale and Cinglani \cite{A-B-C} consider concentration phenomena at isolated local minima and maxima with polynomial degeneracy.

\vskip0.1in

However we remark that the uniqueness and non-degeneracy of the ground state solutions of the limit equation (\ref{q3}) are, in general, rather difficult to prove. They are known so far only for a rather restricted class of nonlinearities $f$. In \cite{Rab}, without the uniqueness and non-degeneracy condition, Rabinowitz proves, by a mountain pass argument, the existence of positive solutions of (\ref{q1}) for small $\e>0$ whenever
$$
\liminf_{|x|\rg\iy}V(x)>\inf_{\RN}V(x).
$$
In \cite{Wang}, with $(V1)$-$(V2)$ Wang proves that these solutions concentrate around the global minimum points of $V$ as $\e\rg 0$. Later, Del Pino and Felmer \cite{Felmer} introduced a penalization approach and proved a localized version of the result by Rabinowitz and Wang. They prove the existence of a single-peak solution which concentrates around the minimum points of $V$ in $O$, providing that the nonlinearity $f$ satisfies $(f_1)$-$(f_2)$, $(H1)$ and the so-called $global Ambrosetti$-$Rabinowitz\ condition$ ($(A$-$R)$ for short): for some $\mu>2, 0<\mu\int_0^tf(s)ds<tf(t), t>0.$ Recently, it has been shown in \cite{Jean,byeon} that $(H1)$ and $(A$-$R)$ are not necessary.

%%%%%%%%%%%%%%%%%%%%%%%%%%%%%%%%%%%%%%%%%%%%%%%%%%%%%%%%%%%%%%%%%%%%%%%%%%%%%%%%%%%%%%%%%%%%%%%%%
%%%%%%%%%%%%%%%%%%%%%%%%%%%%%%%%%%%%%%%%%%%%%%%%%%%%%%%%%%%%%%%%%%%%%%%%%%%%%%%%%%%%%%%%%%%%%%%%%
%%%%%%%%%%%%%%%%%%%%%%%%%%%%%%%%%%%%%%%%%%%%%%%%%%%%%%%%%%%%%%%%%%%%%%%%%%%%%%%%%%%%%%%%%%%%%%%%%

\vskip0.1in

To sum  up,   for the critical case,  the   concentration phenomenon of problem (\ref{q1}) has not been studied so far by variational methods.
Thus similar  to the  subcritical case, when $f$ is critical, it seems natural to expect that there also exists a corresponding solution  to  the  singularly perturbed problem (\ref{q1}) for small $\e>0$ and the similar concentration phenomenon  occurs. In the present paper, we will adopt the ideas of Byeon \cite{byeon} to find solutions of problem (\ref{q1}) in some neighborhood of the set of ground states for problem (\ref{q3}). But it should be stressed that the compactness is the main difficulty in extending the quoted results to critical problems. For critical variational problem, $(PS)$-$condition$ fails. About this aspect, we refer to \cite{Wei,Brezis2,Wang,Alves1,Ruiz} and the references therein. Therefore, the method of Byeon \cite{byeon} is not used directly and some more tricks  are given.

%%%%%%%%%%%%%%%%%%%%%%%%%%%%%%%%%%%%%%%%%%%%%%%%%%%%%%%%%%%%%%%%%%%%%%%%%%%%%%%%%%%%%%%%%%%%%%%%%
%%%%%%%%%%%%%%%%%%%%%%%%%%%%%%%%%%%%%%%%%%%%%%%%%%%%%%%%%%%%%%%%%%%%%%%%%%%%%%%%%%%%%%%%%%%%%%%%%
%%%%%%%%%%%%%%%%%%%%%%%%%%%%%%%%%%%%%%%%%%%%%%%%%%%%%%%%%%%%%%%%%%%%%%%%%%%%%%%%%%%%%%%%%%%%%%%%%
%%%%%%%%%%%%%%%%%%%%%%%%%%%%%%%%%%%%%%%%%%%%%%%%%%%%%%%%%%%%%%%%%%%%%%%%%%%%%%%%%%%%%%%%%%%%%%%%%
%%%%%%%%%%%%%%%%%%%%%%%%%%%%%%%%%%%%%%%%%%%%%%%%%%%%%%%%%%%%%%%%%%%%%%%%%%%%%%%%%%%%%%%%%%%%%%%%%
%%%%%%%%%%%%%%%%%%%%%%%%%%%%%%%%%%%%%%%%%%%%%%%%%%%%%%%%%%%%%%%%%%%%%%%%%%%%%%%%%%%%%%%%%%%%%%%%%

\s{Proof of Theorem $(\ref{Theorem 1})$}

\renewcommand{\theequation}{2.\arabic{equation}}

To study (\ref{q1}), it suffices to study (\ref{q2}).
Let $H_\e$ be the completion of $C_0^\iy(\RN)$ with respect to the norm $$\|u\|_\e=\left(\int_{\RN} |\na u|^2+V_\e u^2\right)^{\frac{1}{2}}.$$ We define a norm $\|\cdot\|$ on $H^1(\RN)$ by $$\|u\|^2=\int_{\RN} |\na u|^2+V_0u^2.$$
Since $\inf_{\RN}V(x)=V_0>0$, we have $H_\e\subset H^1(\RN)$. For any set $B\subset \RN$ and $\e>0$, we define $B_\e\equiv\{x\in\RN: \e x\in B\}$. For $u\in H_\e$, let
$$
P_\e(u)=\frac{1}{2}\int_{\RN} |\na u|^2+V_\e u^2-\int_{\RN} F(u).
$$
Fixing an arbitrary $\mu>0$, we define
\be
\chi_\e(x)=
\bcs\nonumber
0,\ \ \ \ \mbox{if}\ \ x\in O_\e,\\
\varepsilon^{-\mu},\ \ \mbox{if}\ \ x\in \RN\setminus O_\e,
\ecs
\ee
and
$$
Q_\e(u)=\left(\int_{\RN}\chi_\e u^2dx-1\right)_+^2.
$$
The functional $Q_\e$ will act as a penalization to force the concentration phenomena to occur inside $O$. This type of penalization was first introduced in $\cite{byeon1}$. Finally, let $\G_\e:H_\e\rg \R$ be given by $$\G_\e(u)=P_\e(u)+Q_\e(u).$$ It is standard to show that $\G_\e\in C^1(H_\e)$. From now on we may assume that $f(t) = 0$ for $t \le 0.$ In this case any critical point of $\G$ is positive by the maximum principle. Clearly a critical point of $P_\e$ corresponds to a solution of $(\ref{q2})$. To find solutions of $(\ref{q2})$ which concentrate in $O$ as $\e\rg 0$, we shall search critical points of $\G_\e$ for which $Q_\e$ is zero.

First, we study some properties of the solutions of $(\ref{q3})$. Without loss of generality, we may assume that $0\in \mathcal{M}$. For any set $B\subset \RN$ and $\dd>0$, we define $B^{\dd}\equiv\{x\in\RN|dist(x,B)\le\dd\}$. As we already mentioned, the following equations for $a>0$ are limiting equations of $(\ref{q2})$
\be\lab{q4}
-\DD u+au=f(u),\ \ u>0,\ \ u\in H^1(\RN).
\ee
We define an energy functional for the limiting problems $(\ref{q4})$ by
$$
L_a(u)=\frac{1}{2}\int_{\RN}|\na u|^2+au^2-\int_{\RN}F(u),\ \ u\in H^1(\RN),
$$
where $F(s)=\int_0^sf(t)dt$. In $\cite{Zhang-Zou}$, we proved that, if $p>2,N\ge 4$ or $p>4, N=3$ and $(F1)$-$(F3)$ hold, there exists a least energy solution of $(\ref{q4})$ for any $a>0$. Moreover, each such solution $U$ of $(\ref{q4})$ satisfies Pohozaev's identity
\be\label{eq}
\frac{N-2}{2}\int_{\RN}|\na U|^2\,dx=N\int_{\RN}\left(F(U)-\frac{a}{2}U^2\right)\,dx,
\ee
and so \be\label{Poho}\int_{\RN}|
\nabla U|^2\,dx=N L_a(U).\ee

Let $S_a$ be the set of least energy solutions $U$ of $(\ref{q4})$ satisfying $U(0)=\max_{x\in \RN}U(x)$.  The following result on $S_a$ was proved in \cite{Byeon-Zhang-Zou}.
\bo(see \cite[Proposition 2.1]{Byeon-Zhang-Zou})\lab{prop2.1}\
\begin{itemize}
\item [$(1)$] For any $U\in S_a$, $U$ is radially symmetric and $\frac{\partial U}{\partial r}<0$ for all $r>0$.
\item [$(2)$] $S_a$ is compact in $H^1(\RN)$.
\item [$(3)$] $0<\inf\{\|U\|_{\infty}:U\in S_a\}\le \sup\{\|U\|_{\infty}:U\in S_a\}<\iy$.
\item [$(4)$] There exist $C,c>0$, independent of $U\in S_a$, such that $|D^{\al}U(x)|\le C\exp(-c|x|), \,x\in \RN$ for $|\al|=0,1$.
\end{itemize}
\eo

Let $E_m=L_m(U)$ for $U\in S_m$ and $10\delta=dist(\mathcal{M},O^c)$. We fix a $\beta\in (0,\delta)$ and a cut-off $\vp\in C_0^\infty(\RN)$ such that $0\le\vp\le1,\vp(x)=1$ for $|x|\le \beta$ and $\vp(x)=0$ for $|x|\ge 2\beta$. Let $\vp_{\e}(y)=\vp(\e y), y\in\RN$ and for each $x\in \mathcal{M}^{\beta}$ and $U\in S_m$, we define
$$
U_{\e}^x(y):=\vp_{\e}\left(y-\frac{x}{\e}\right)U\left(y-\frac{x}{\e}\right).
$$
We will find a solution near the set
$$
X_{\e}:=\{U_{\e}^x(y)\,\,|\,\,x\in \mathcal{M}^{\beta}, U\in S_m\}.
$$
for sufficiently small $\e>0$. We note that $0\in \mathcal{M}$ and define
$$
W_{\e}(y)=\vp_{\e}(y)U(y),
$$
where $U\in S_m$ is arbitrary but fixed. Setting $W_{\e,t}(y)=\vp_{\e}(y)U(\frac{y}{t})$, we see that $\Gamma_{\e}(W_{\e,t})=P_{\e}(W_{\e,t})$ for $t\ge 0$. By (\ref{eq}), for $U_t(x)=U(\frac{x}{t})$ we have
$$
L_m(U_t)=\left(\frac{t^{N-2}}{2}-\frac{(N-2)t^N}{2N}\right)\int_{\RN}|\nabla U|^2.
$$
Thus, there exists $t_0>1$ such that $L_m(U_t)<-2$ for $t\ge t_0$.

Finally, we define a min-max value $C_\e$:
$$C_{\e}:=\inf_{\gamma\in \Phi_{\e}}\max_{s\in [0,1]}\Gamma_{\e}(\gamma(s)),$$
where $\Phi_{\e}:=\{\gamma \in C([0,1],H_{\e}): \gamma(0)=0,\gamma(1)=W_{\e,t_0}\}$. We can check that $\Gamma_{\e}(\gamma(1))<-2$ for any $\e>0$ sufficiently small. Let $\g_\e(s)=W_{\e,st_0}$ for $s\in (0,1]$ and $\g_\e(0)=0$. We denote
$$
D_\e:=\max_{s\in [0,1]}\G_\e(\g_\e(s)).
$$
Recalling that in \cite{Zhang-Zou}, the authors proved that, for equation (\ref{q3}) the mountain pass level corresponds to the least energy level. Then similar as in \cite{byeon}, we can prove that

\bo\lab{prop2.2}
$\lim\limits_{\e\rg 0}C_{\e}=\lim\limits_{\e\rg 0}D_{\e}=E_m.$
\eo
Now define
$$
\G_\e^\al:=\{u\in H_\e: \G_\e(u)\le\al\}
$$
and for a set $A\subset H_\e$ and $\al>0$, let
$$
A^\al:=\{u\in H_\e: \inf_{v\in A}\|u-v\|_\e\le\al\}.
$$
The following lemma was introduced in \cite{Wang} and will be used in the following proof.
\bl(\cite{Wang})
\lab{Wang}
Let $R$ be a positive number and $\{u_n\}$ a bounded sequence in $H^1(\RN)$. If
$$
\lim_{n\rg \iy}\sup_{x\in\RN}\int_{B(x,R)}|u_n|^{2^{\ast}}=0,
$$
then $u_n\rg 0$ in $L^{2^{\ast}}(\RN)$ as $n\rg \iy$.
\el

To continue our proof, we need the following lemma from Benci and Cerami \cite{Benci}.

\bl(see \cite[Lemma 2.7]{Benci})
\lab{Benci}
Let $\{u_m\}\subset H^1_{loc}(\RN)$ be a sequence of functions such that
$$u_m\rightharpoonup 0 \,\,\,\hbox{weakly in $H^1(\RN)$}.$$
Suppose that there exist a bounded open set $Q\subset \RN$ and a positive constant $\gamma>0$ such that
$$
\int_{Q}|\nabla u_m|^2\,dx\ge \gamma>0, \quad \int_{Q}|u_m|^{2^\ast}\ge \gamma>0.
$$
Moreover suppose that
$$\Delta u_m+ |u_m|^{2^\ast-2} u_m=\chi_m,$$
where $\chi_m\in H^{-1}(\RN)$ and
$$|\langle\chi_m, \phi\rangle|\le \e_m\|\phi\|_{H^1(\RN)},\quad\forall\,\phi\in C^{\iy}_0(U),$$
where $U$ is an open neighborhood of $Q$ and $\e_m$ is a sequence converging to $0$. Then there exist a sequence of points $\{y_m\}\subset\RN$ and
a sequence of positive numbers $\{\sigma_m\}$ such that
$$v_m(x):=\sigma_m^{(N-2)/2}u_m(\sigma_m x+y_m)$$
converges weakly in $D^{1,2}(\RN)$ to a nontrivial solution $v$ of
$$-\Delta u=|u|^{2^\ast-2}u,\quad u\in D^{1,2}(\RN).$$
Moreover,
$$y_m\to \bar{y}\in\overline{Q}\quad\hbox{and}\quad \sigma_m\to 0.$$
\el

The following proposition is very important, and its proof is much more delicate in case of critical growth.
\bo\lab{prop2.4}
Let $\{\e_i\}_{i=1}^\infty$ be such that $\lim_{i\rg \infty}\e_i=0$ and $\{u_{\e_i}\}\subset X_{\e_i}^d$ such that
$$
\lim_{i\rg\infty}\Gamma_{\e_i}(u_{\e_i})\le E_m\ \mbox{and}\ \lim_{i\rg\infty}\Gamma_{\e_i}^{'}(u_{\e_i})=0.
$$
Then for sufficiently small $d>0$, there exits, up to a subsequence, $\{y_i\}_{i=1}^\infty\subset\RN, x\in \mathcal{M}, U\in S_m$ such that
$$
\lim_{i\rg \infty}|\e_iy_i-x|=0\ \mbox{and}\ \lim_{i\rg\infty}\|u_{\e_i}-\varphi_{\e_i}(\cdot-y_i)U(\cdot-y_i)\|_{\e_i}=0.
$$
\eo
\bp For convenience, we write $\e$ for $\e_i$. By the definition of $X_\e^d$, there exist $\{U_\e\}\subset S_m$ and $\{x_\e\}\subset\mathcal{M}^\bb$ with
$$
\|u_\e-\vp_\e(\cdot-\frac{x_\e}{\e})U_\e(\cdot-\frac{x_\e}{\e})\|_\e\le d.
$$
Since $S_m$ and $\mathcal{M}^\bb$ are compact, there exist $Z\in S_m, x\in \mathcal{M}^\bb$ such that $U_\e\rg Z$ in $H^1(\RN)$ and $x_\e\rg x$. Thus, for small $\e>0$,
\be\lab{ff1}
\|u_\e-\vp_\e(\cdot-\frac{x_\e}{\e})Z(\cdot-\frac{x_\e}{\e})\|_\e\le 2d.
\ee

\noindent {\bf Step 1.} We claim that
$$
\liminf_{\e\rg 0}\sup_{y\in A_\e}\int_{B(y,1)}|u_\e|^{2^{\ast}}=0,
$$
where $A_\e=B(\frac{x_\e}{\e},\frac{3\bb}{\e})\setminus B(\frac{x_\e}{\e},\frac{\bb}{2\e})$.
If the claim is true, by Lemma \ref{Wang} we see that
$$
u_\e\rg 0\ \ \mbox{in}\ \  L^{2^{\ast}}(B_\e),
$$
where $B_\e=B(\frac{x_\e}{\e},\frac{2\bb}{\e})\setminus B(\frac{x_\e}{\e},\frac{\bb}{\e})$.

By Lemma \ref{Wang}, assume by contradiction that there exists $r>0$, such that
$$
\liminf_{\e\rg 0}\sup_{y\in A_\e}\int_{B(y,1)}|u_\e|^{2^{\ast}}=2r>0,
$$
then there exists $y_\e\in A_\e$ such that for small $\e>0$, $\int_{B(y_\e,1)}|u_\e|^{2^{\ast}}\ge r$.
Note that $y_\e\in A_\e$ and there exists $x_0\in \mathcal{M}^{4\bb}\subset O$ such that $\e y_\e\rg x_0$. Let $v_\e(y)=u_\e(y+y_\e)$, then, for $\e$ small,
\be\lab{ff2}
\int_{B(0,1)}|v_\e|^{2^{\ast}}\ge r
\ee
and up to a subsequence, $v_\e\rg v$ weakly in $H^1(\RN)$ and $v$ satisfies
$$
-\DD v+V(x_0)v=f(v)\ \ \mbox{in}\ \ \RN.
$$

\vskip0.1in

{\bf Case 1.} If $v\not\equiv0$, then for sufficiently large $R>0$,
$$
\liminf_{\e\rg 0}\int_{B(y_\e,R)}|\na u_\e|^2\ge \frac{1}{2}\int_{\RN}|\na v|^2=\frac{N}{2}L_{V(x_0)}(v).
$$
By definition, $L_{V(x_0)}(v)\ge E_{V(x_0)}$. Now recalling from \cite{{Zhang-Zou}} that $E_a>E_b$ if $a>b$. Then since $x_0\in O$, we have $V(x_0)\ge m$ and $\liminf_{\e\rg 0}\int_{B(y_\e,R)}|\na u_\e|^2\ge \frac{N}{2}E_m>0$, which is a contradiction with (\ref{ff1}) if $d$ is small enough.
\vskip0.1in

{\bf Case 2.} If $v\equiv0$, i.e., $v_\e\rightharpoonup 0$ weakly in $H^1(\RN)$, then $v_\e\rg 0$ strongly in $L_{loc}^q(\RN)$ for $q\in [2,2^{\ast})$. Thus, by (\ref{ff2}) and Sobolev' embedding theorem, there exists $C>0$ (independent of $\e$) such that, for $\e$ small,
\be\lab{ff3}
\int_{B(0,1)}|\na v_\e|^2\ge Cr^{\frac{2}{2^{\ast}}}>0.
\ee
Now, we claim that
\be\lab{ff4}
\lim_{\e\rg 0}\sup_{\stackrel{\phi\in C_0^\iy(\Omega)}{\|\phi\|=1}}|\langle\rho_\e,\phi\rangle|=0,
\ee
where $\Omega=B(0,2), \rho_\e\in H^{-1}(\RN)$ and $\rho_\e=\DD v_\e+|v_\e|^{2^{\ast}-2}v_\e$.
For $\e>0$ small enough, it is easy to check that $\int_{\RN}\chi_\e u_\e\phi(\cdot-y_\e)\equiv0$ uniformly for any $\phi\in C_0^\iy(\Omega)$. Thus, for any $\phi\in C_0^\iy(\Omega)$ and $\|\phi\|=1$,
{\allowdisplaybreaks
\begin{align}
\lab{ff5}
\langle\rho_\e,\phi\rangle=&\int_{\RN}V_\e u_\e \phi(\cdot-y_\e)-\langle\G_\e'(u_\e),\phi(\cdot-y_\e)\rangle\nonumber\\
&-\int_{\RN}\big(f(u_\e)-|u_\e|^{2^{\ast}-2}u_\e\big)\phi(\cdot-y_\e)\nonumber\\
:=&J_1+J_2+J_3.
\end{align}
}%
First, since $\sup_{y\in \Omega}V_\e(y+y_\e)<\iy$ uniformly for small $\e$ and $v_\e\rg 0$ strongly in $L_{loc}^2(\RN)$, we have
\begin{align*}
|J_1|&\le \sup_{y\in \Omega}V_\e(y+y_\e)\left(\int_{\Omega}|v_\e|^2\right)^{\frac{1}{2}}\left(\int_{\Omega}|\phi|^2\right)^{\frac{1}{2}}\\
&\to 0 \quad\hbox{as $\e\to 0$, uniformly for $\phi\in C^{\iy}_0(\Omega), \|\phi\|=1$.}
\end{align*}
Second, by $(F1)$-$(F2)$, we see that $\lim_{|s|\rg 0}g(s)/|s|=\lim_{|s|\rg \iy}g(s)/|s|^{2^{\ast}-1}=0$, where $g(s)=f(s)-|s|^{2^{\ast}-2}s$. Let $\ti{g}(s)=|g(s)|^{\frac{2^{\ast}}{2^{\ast}-1}}$, then
$$
\lim_{s\rg 0}\frac{\ti{g}(s)}{|s|^{\frac{2^{\ast}}{2^{\ast}-1}}}=0\ \ \ \mbox{and}\ \ \ \ \lim_{s\rg \iy}\frac{\ti{g}(s)}{|s|^{2^{\ast}}}=0.
$$
Since $v_\e\rightharpoonup 0$ weakly in $H^1(\RN)$, $\sup_{\e>0}\int_{\Omega} |v_\e|^{\frac{2^{\ast}}{2^{\ast}-1}}+|v_\e|^{2^{\ast}}<\iy$. Then, it follows from $v_\e\rg 0$ strongly in $L_{loc}^2(\RN)$ and the compactness lemma of Strauss \cite{Strauss} that $\ti{g}(v_\e)\rg 0$ strongly in $L^1(\Omega)$ as $\e\rg 0$. Thus, by Lebesgue dominated convergence theorem,
$$
\lim_{\e\rg 0}|J_3|\le\lim_{\e\rg 0}\left(\int_{\Omega}\ti{g}(v_\e)\right)^{\frac{2^{\ast}-1}{2^{\ast}}}
\left(\int_{\Omega}|\phi|^{2^{\ast}}\right)^{\frac{1}{2^{\ast}}}=0
$$
uniformly for $\phi\in C^{\iy}_0(\Omega), \|\phi\|=1$.
Combining these with $\G_\e'(u_\e)\rg 0$, we see that the claim (\ref{ff4}) is true. By Lemma \ref{Benci}, we see from (\ref{ff2})-(\ref{ff4}) that, there exists $\ti{y}_\e\in\RN$ and $\sigma_\e>0$, such that $\ti{y}_\e\rg \ti{y}\in \overline{B(0,1)}$, $\sigma_\e\rg 0$
and
$$
w_\e(y):=\sigma_\e^{\frac{N-2}{2}}v_\e(\sigma_\e y+\ti{y}_\e)
$$
converges weakly to $w$ in $D^{1,2}(\RN)$, where $w$ is a nontrivial solution of $-\DD v=|v|^{2^{\ast}-2}v$ in $D^{1,2}(\RN)$. Then there exists $R>0$, such that $\int_{B(0,R)}|\na w|^2\ge\frac{S^{\frac{N}{2}}}{2}$. Hence, for $\e>0$ small enough,
\begin{eqnarray*}
\liminf_{\e\rg 0}\int_{B(y_\e,2)}|\na u_\e|^2\ge\liminf_{\e\rg 0}\int_{B(\ti{y}_\e,\sigma_\e R)}|\na v_\e|^2\ge\int_{B(0,R)}|\na w|^2\ge\frac{S^{\frac{N}{2}}}{2},
\end{eqnarray*}
which is a contradiction with (\ref{ff1}) if $d>0$ is small enough. Therefore, Step 1 is proved.

\vskip0.1in

\noindent {\bf Step 2.} Let $u_\e^1(y)=\varphi_\e(y-\frac{x_\e}{\e})u_\e(y), u_\e^2=u_\e-u_\e^1$. We claim that, for small $d>0$, $\G_\e(u_\e^2)\ge 0$ and
$$
\G_\e(u_\e)\ge \G_\e(u_\e^1)+\G_\e(u_\e^2)+o(1),\ \ \mbox{as}\ \ \e\rg 0.
$$
Similarly as in \cite[Proposition 4]{byeon}, for small $d>0$, $\G_\e(u_\e^2)\ge 0$. By $(F1)$, for any $\delta>0$, there exits $C_\delta>0$ such that $|F(s)|\le \delta|s|^2+C_
\delta|s|^{2^{\ast}}$ for $s\in\R$. Then by $u_\e\rg 0$ strongly in $L^{2^{\ast}}(B_\e)$ which has been proved in Step 1, we have
{\allowdisplaybreaks
\begin{align*}
\limsup_{\e\rg 0}\left|\int_{\RN}F(u_\e)-F(u_\e^1)-F(u_\e^2)\right|
=&\limsup_{\e\rg 0}\left|\int_{B_\e}F(u_\e)-F(u_\e^1)-F(u_\e^2)\right|\\
\le&\limsup_{\e\rg 0}\int_{B_\e}\delta|u_\e|^2+C_
\delta|u_\e|^{2^{\ast}}\\
\le&C\delta,
\end{align*}
}%
where $C$ is a positive constant (independent of $\e, \delta$). Since $\delta$ is arbitrary, $\int_{\RN}F(u_\e)-F(u_\e^1)-F(u_\e^2)=o(1)$ as $\e\rg 0$. Then the claim can be proved similarly as in \cite[Propostion 4]{byeon}. We omit the details.

\vskip0.1in

\noindent {\bf Step 3.} Let $w_\e(y):=u^1_\e(y+\frac{x_\e}{\e})=\varphi_\e(y)u_\e(y+\frac{x_\e}{\e})$. Up to a subsequence, $w_\e\rightharpoonup w$ weakly in $H^1(\RN)$, $w_\e\rg w$ a.e. in $\RN$. Now, we claim that
$$
w_\e\rg w\ \ \mbox{strongly in}\ \ L^{2^{\ast}}(\RN).
$$
By Lemma \ref{Wang}, assume by contradiction that there exists $r>0$, such that
$$
\liminf_{\e\rg 0}\sup_{z\in \RN}\int_{B(z,1)}|w_\e-w|^{2^{\ast}}=2r>0.
$$
Then, there exists $z_\e\in\RN$ such that $\liminf_{\e\rg 0}\int_{B(z_\e,1)}|w_\e-w|^{2^{\ast}}>r$.

\vskip0.1in

{\bf Case 1.} $\{z_\e\}$ is bounded, i.e., $|z_\e|\le a$ for some $a>0$. Then for $\e$ small,
\be\lab{gg1}
\int_{B(0,a+1)}|v_\e|^{2^{\ast}}>r,
\ee
where $v_\e=w_\e-w$ and $v_\e\rightharpoonup 0$ weakly in $H^1(\RN)$. Similar as in Step 1, there exists $C>0$ (independent of $\e$) such that, for $\e$ small,
\be\lab{gg2}
\int_{B(0,a+1)}|\na v_\e|^2\ge Cr^{\frac{2}{2^{\ast}}}>0.
\ee
Now, we claim that
\be\lab{gg3}
\lim_{\e\rg 0}\sup_{\stackrel{\phi\in C_0^\iy(\Omega)}{\|\phi\|=1}}|\langle\rho_\e,\phi\rangle|=0,
\ee
where $\Omega=B(0,a+2), \rho_\e\in H^{-1}(\RN)$ and $\rho_\e=\DD v_\e+|v_\e|^{2^{\ast}-2}v_\e$.
For $\e>0$ small enough, it is easy to check that $\int_{\RN}\chi_\e u_\e\phi(\cdot-\frac{x_\e}{\e})\equiv0$ uniformly for any $\phi\in C_0^\iy(\Omega)$. Thus, we see that for $\e$ small and any $\phi\in C_0^\iy(\Omega)$, $$\langle\G_\e'(u_\e),\phi(\cdot-\frac{x_\e}{\e})\rangle=\int_{\RN}\na w_\e\na \phi+V_\e(y+\frac{x_\e}{\e})w_\e\phi-f(w_\e)\phi.$$
Then,
\be\lab{gg4}
\int_{\RN}\na w_\e\na \phi+V_\e(y+\frac{x_\e}{\e})w_\e\phi-f(w_\e)\phi=o(1),
\ee
uniformly for $\phi\in C_0^\iy(\Omega), \|\phi\|=1$. Noting that $w_\e\rightharpoonup w$ weakly in $H^1(\RN)$ and $x_\e\rg x$, by a standard way we can see that $w$ satisfies $-\DD w(y)+V(x)w(y)=f(w(y))$ in $\RN$. Thus,
\be\lab{gg5}
\int_{\RN}\na w\na \phi+V(x)w\phi-f(w)\phi=0,\ \ \mbox{for}\ \ \phi\in C_0^\iy(\Omega).
\ee
Moreover, by elliptic estimates $w\in L^{\iy}(\RN)$. Then it follows that
\be\lab{gg6}
\int_{\RN}|w_\e|^{2^{\ast}-2}w_\e\phi-|v_\e|^{2^{\ast}-2}v_\e\phi-|w|^{2^{\ast}-2}w\phi=o(1),
\ee
uniformly for $\phi\in C_0^\iy(\Omega), \|\phi\|=1$.
It follows from (\ref{gg4})-(\ref{gg6}) that
{\allowdisplaybreaks
\begin{align}
\lab{gg7}
\langle\rho_\e,\phi\rangle&=\int_{\Omega}\big(V(\e y+x_\e)w_\e-V(x)w\big)\phi
+\int_{\Omega}\big(g(w)-g(w_\e)\big)\phi+o(1)\nonumber\\
&=:K_1+K_2,
\end{align}
}%
where $g(s)=f(s)-|s|^{2^{\ast}-2}s$ and $o(1)\rg 0$ as $\e\rg 0$ uniformly for $\phi\in C_0^\iy(\Omega), \|\phi\|=1$.
Noting that $w_\e\rg w$ strongly in $L^2(\Omega)$ and $x_\e\rg x$, it is easy to check that $K_1=o(1)$ as $\e\rg 0$ uniformly for $\phi\in C_0^\iy(\Omega), \|\phi\|=1$. In the following, we will show that
\be\lab{gg8}
\lim_{\e\rg 0}\int_{\Omega}|g(w_\e)-g(w)|^{\frac{2^{\ast}}{2^{\ast}-1}}=0.
\ee
Since $w_\e\rg w$ weakly in $H^1(\RN)$, $\sup_{\e>0}\int_{\Omega} |w_\e|^{\frac{2^{\ast}}{2^{\ast}-1}}+|w_\e|^{2^{\ast}}<\iy$, similar as in Step 1, $\ti{g}(w_\e)\rg \ti{g}(w)$ strongly in $L^1(\Omega)$ as $\e\rg 0$. Thus,  by Lebesgue dominated convergence theorem, (\ref{gg8}) is proved. Then,
$$
\lim_{\e\rg 0}|K_2|\le\lim_{\e\rg 0}\left(\int_{\Omega}|g(w_\e)-g(w)|^{\frac{2^{\ast}}{2^{\ast}-1}}\right)^{\frac{2^{\ast}-1}{2^{\ast}}}
\left(\int_{\Omega}|\phi|^{2^{\ast}}\right)^{\frac{1}{2^{\ast}}}=0.
$$
Therefore, (\ref{gg3}) follows from (\ref{gg7}). By Lemma \ref{Benci} again, we see from (\ref{gg1})-(\ref{gg3}) that, there exist $\ti{z}_\e\in\RN$ and $\sigma_\e>0$, such that $\ti{z}_\e\rg \ti{z}\in \overline{B(0,a+1)}$, $\sigma_\e\rg 0$
and
$$
\ti{w}_\e(y):=\sigma_\e^{\frac{N-2}{2}}v_\e(\sigma_\e y+\ti{z}_\e)
$$
converges weakly to $\ti{w}$ in $D^{1,2}(\RN)$, where $\ti{w}$ is a nontrivial solution of $-\DD \ti{w}=|\ti{w}|^{2^{\ast}-2}\ti{w}$ in $D^{1,2}(\RN)$. Then,
{\allowdisplaybreaks
\begin{align*}
\int_{\RN}|\na \ti{w}|^2&\le\liminf_{\e\rg 0}\int_{\RN}|\na \ti{w}_\e|^2
=\liminf_{\e\rg 0}\int_{\RN}|\na v_\e|^2\\
&=\liminf_{\e\rg 0}\int_{\RN}|\na w_\e|^2-|\na w|^2.
\end{align*}
}%
Noting that $\int_{\RN}|\na \ti{w}|^2\ge S^{\frac{N}{2}}$, we see that $\liminf_{\e\rg 0}\int_{\RN}|\na w_\e|^2\ge S^{\frac{N}{2}}$. It follows that $\liminf_{\e\rg 0}\int_{\RN}|\na u_\e|^2\ge S^{\frac{N}{2}}$. By (\ref{ff1}), we have
{\allowdisplaybreaks
\begin{align*}
\left(\int_{\RN}|\na Z|^2\right)^{\frac{1}{2}}&=\lim_{\e\rg 0}\left(\int_{\RN}\left|\na (\varphi_\e Z)\right|^2\right)^{\frac{1}{2}}
\ge\liminf_{\e\rg 0}\left(\int_{\RN}|\na u_\e|^2\right)^{\frac{1}{2}}-2d\\
&\ge\left(S^{\frac{N}{2}}\right)^{\frac{1}{2}}-2d.
\end{align*}
}%
Recalling that $E_m<\frac{S^{\frac{N}{2}}}{N}$, we see that $\int_{\RN}|\na Z|^2>N E_m$ for small $d>0$, which is a contradiction, since $Z\in S_m$.

\vskip0.1in

{\bf Case 2.} $\{z_\e\}$ is unbounded. Without loss of generality, $\lim_{\e\rg 0}|z_\e|=\iy$. Then, $\liminf_{\e\rg 0}\int_{B(z_\e,1)}|w_\e|^{2^{\ast}}\ge r$, i.e.,
$$\liminf_{\e\rg 0}\int_{B(z_\e,1)}|\varphi_\e(y)u_\e(y+\frac{x_\e}{\e})|^{2^{\ast}}\ge r.$$
Since $\varphi(y)=0$ for $|y|\ge 2\bb$, we see that $|z_\e|\le \frac{3\bb}{\e}$ for $\e$ small. If $|z_\e|\ge \frac{\bb}{2\e}$ for $\e$ small, then $z_\e\in B(0,\frac{3\bb}{\e})\setminus B(0,\frac{\bb}{2\e})$ and by Step 1, we get that
{\allowdisplaybreaks
\begin{align*}
\liminf_{\e\rg 0}\int_{B(z_\e,1)}|w_\e|^{2^{\ast}}
\le&\liminf_{\e\rg 0}\sup_{z\in B(0,\frac{3\bb}{\e})\setminus B(0,\frac{\bb}{2\e})}\int_{B(z,1)}|u_\e(y+\frac{x_\e}{\e})|^{2^{\ast}}\\
=&\liminf_{\e\rg 0}\sup_{z\in A_\e}\int_{B(z,1)}|u_\e|^{2^{\ast}}\\
=&0,
\end{align*}
}%
which is a contradiction. Thus, $|z_\e|\le \frac{\bb}{2\e}$ for $\e$ small. Assume that $\e z_\e\rg z_0\in \overline{B(0,\frac{\bb}{2})}$ and $\ti{w}_\e=w_\e(y+z_\e)\rightharpoonup \ti{w}$ weakly in $H^1(\RN)$. If $\ti{w}\not\equiv 0$, we can see that $\ti{w}$ satisfies
$$
-\DD \ti{w}(y)+V(x+z_0)\ti{w}(y)=f(\ti{w}(y))\ \ \mbox{in}\ \  \RN.
$$
Similar as in Step 1, we get a contradiction if $d>0$ is small enough. Thus, $\ti{w}\equiv 0$, i.e., $\ti{w}_\e\rg 0$ weakly in $H^1(\RN)$. Meanwhile,
\be\lab{kk2}
\int_{B(0,1)}|\ti{w}_\e|^{2^{\ast}}\ge r>0
\ee
and there exists $C>0$ (independent of $\e$) such that, for $\e$ small,
\be\lab{kk3}
\int_{B(0,1)}|\na \ti{w}_\e|^2\ge Cr^{\frac{2}{2^{\ast}}}>0.
\ee
Now, we claim that
\be\lab{kk4}
\lim_{\e\rg 0}\sup_{\stackrel{\phi\in C_0^\iy(\Omega)}{\|\phi\|=1}}|\langle\ti{\rho}_\e,\phi\rangle|=0,
\ee
where $\Omega=B(0,2), \ti{\rho}_\e\in H^{-1}(\RN)$ and $\ti{\rho}_\e=\DD \ti{w}_\e+|\ti{w}_\e|^{2^{\ast}-2}\ti{w}_\e$.
For $\e>0$ small enough, it is easy to check that $\int_{\RN}\chi_\e u_\e\phi(\cdot-z_\e-\frac{x_\e}{\e})\equiv0$ and $\ti{w}_\e(y)=u_\e(y+z_\e+\frac{x_\e}{\e})$ uniformly for any $y\in \Omega$ and $\phi\in C_0^\iy(\Omega)$. Thus, for any $\phi\in C_0^\iy(\Omega)$ and $\|\phi\|=1$,
\begin{align*}
\langle\ti{\rho}_\e,\phi\rangle=&\int_{\RN}V_\e u_\e \phi(\cdot-z_\e-\frac{x_\e}{\e})-\langle\G_\e'(u_\e),\phi(\cdot-z_\e-\frac{x_\e}{\e})\rangle\\
&-\int_{\RN}\big(f(u_\e)-|u_\e|^{2^{\ast}-2}u_\e\big)\phi(\cdot-z_\e-\frac{x_\e}{\e}).
\end{align*}
Similar as in Step 1, (\ref{kk4}) is true. By Lemma \ref{Benci} again, we see from (\ref{kk2})-(\ref{kk4}) that, there exist $\ti{y}_\e\in\RN$ and $\sigma_\e>0$, such that $\ti{y}_\e\rg \ti{y}\in \overline{B(0,1)}$, $\sigma_\e\rg 0$
and
$$
\widehat{w}_\e(y):=\sigma_\e^{\frac{N-2}{2}}\ti{w}_\e(\sigma_\e y+\ti{y}_\e)
$$
converges weakly to $\widehat{w}$ in $D^{1,2}(\RN)$, where $\widehat{w}$ is a nontrivial solution of $-\DD v=|v|^{2^{\ast}-2}v$ in $D^{1,2}(\RN)$. Note that there exists $R>0$, such that $\int_{B(0,R)}|\na \widehat{w}|^2\ge\frac{S^{\frac{N}{2}}}{2}$. Then, for $\e>0$ small enough,
\begin{eqnarray*}
\liminf_{\e\rg 0}\int_{B(z_\e+\frac{x_\e}{\e},2)}|\na u_\e|^2\ge\liminf_{\e\rg 0}\int_{B(\ti{y}_\e,\sigma_\e R)}|\na \ti{w}_\e|^2\ge\int_{B(0,R)}|\na \widehat{w}|^2\ge\frac{S^{\frac{N}{2}}}{2},
\end{eqnarray*}
which is a contradiction with (\ref{ff1}) if $d>0$ is small enough. Therefore, $w_\e\rg w\ \ \mbox{strongly in}\ \ L^{2^{\ast}}(\RN)$.

\vskip0.1in

\noindent {\bf Step 4.} By Step 3, we deduce that
 $$\lim_{\e\to0}\int_{\RN}F(w_\e)\,dx=\int_{\RN}F(w)\,dx.$$
Then similarly as in \cite[Proposition 4]{byeon}, there exist $U\in S_m$ and $y_\e\in \RN$, such that $\lim_{\e\rg 0}|\e y_\e-x|=0$ and $\lim_{\e\rg 0}\|u_\e-\varphi_\e(\cdot-y_\e)U(\cdot-y_\e)\|_\e=0$. This completes the proof.
\ep

\bo\lab{prop2.5}
For sufficiently small $\e>0$ and sufficiently large $R>0$, there exists a sequence $\{u_{n,\e}^R\}_{n=1}^\infty\subset X_\e^d\cap H_0^1(B(0,\frac{R}{\e}))\cap \G_\e^{D_\e}$ such that $|\G_\e'(u_{n,\e}^R)|\rg 0$ in $\big(H_0^1(B(0,\frac{R}{\e}))\big)^{\ast}$.
\eo
\bp
By Proposition \ref{prop2.4}, the proof can be done similarly as in \cite{byeon,byeon2} and the details are omitted here.
\ep
\bo\lab{prop2.6}
For sufficiently small $\e,d>0$, $\G_\e$ has a nontrivial critical point $u_\e\in X_\e^d\cap\G_\e^{D_\e}$.
\eo
\bp
Let $\e>0$ be fixed, small enough and $$d\in \left(0,\frac{1}{2}(\frac{N}{N-2})^{\frac{N-2}{4}}S^{\frac{N}{4}}\right).$$
\vskip0.1in

\noindent {\bf Step 1.} For sufficiently large $R>0$, we claim that $\G_\e$ has a nontrivial critical point $u_\e^R\in X_\e^d\cap H_0^1(B(0,\frac{R}{\e}))\cap \G_\e^{D_\e}$. Proposition \ref{prop2.5} implies that for some $R_0>0$ and any $R>R_0$, there exists a sequence  $\{u_{n,\e}^R\}_{n=1}^\infty\subset X_\e^d\cap H_0^1(B(0,\frac{R}{\e}))\cap \G_\e^{D_\e}$ such that  $|\G_\e'(u_{n,\e}^R)|\rg 0$ in $\big(H_0^1(B(0,\frac{R}{\e}))\big)^{\ast}$. Since $X_\e^d$ is bounded, we can assume that $u_{n,\e}^R\rightharpoonup u_\e^R$ weakly in $H_0^1(B(0,\frac{R}{\e}))$ as $n\rg\infty$. Then it can be proved by a standard way that $\G_\e'(u_\e^R)=0$ in $H_0^1(B(0,\frac{R}{\e}))$. We write that $u_{n,\e}^R=v_{n,\e}^R+w_{n,\e}^R$ with $v_{n,\e}^R\in X_\e$ and $\|w_{n,\e}^R\|_\e\le d$. Since $S_m$ is compact, taking a subsequence if it is necessary, we can assume that $v_{n,\e}^R\rg v_\e^R\in X_\e$ strongly in $H_0^1(B(0,\frac{R}{\e}))$ and $w_{n,\e}^R\rg w_\e^R$ weakly in $H_0^1(B(0,\frac{R}{\e}))$ as $n\rg \infty$. Then we have that $u_\e^R= v_\e^R+w_\e^R$ with $\|w_\e^R\|_\e\le d$, i.e., $u_\e^R\in X_\e^d$. Now, we will prove that $u_\e^R\in \G_\e^{D_\e}$. From the weak convergence of $u_{n,\e}^R$ to $u_\e^R$ in $H_0^1(B(0,\frac{R}{\e}))$ and Brezis-Lieb lemma, we see that as $n\rg \infty$,
\begin{eqnarray*}
&\ &\|\na u_{n,\e}^R\|_2^2=\|\na (u_{n,\e}^R-u_\e^R)\|_2^2+\|\na u_\e^R\|_2^2+o(1),\\
&\ &\|u_{n,\e}^R\|_{2^{\ast}}^{2^{\ast}}=\|u_{n,\e}^R-u_\e^R\|_{2^{\ast}}^{2^{\ast}}+\|u_\e^R\|_{2^{\ast}}^{2^{\ast}}+o(1).
\end{eqnarray*}
Then, by $\lim_{s\rg\infty}\frac{f(s)-s^{2^{\ast}-1}}{s^{2^{\ast}-1}}=0$, we get that as $n\rg \infty$,
$$
\int_{\RN} F(u_{n,\e}^R)-F(u_{n,\e}^R-u_\e^R)-F(u_\e^R)=o(1).
$$
Meanwhile, by the compactness embedding of $H_0^1(B(0,\frac{R}{\e}))\hookrightarrow L^q(B(0,\frac{R}{\e}))(q\in[2,2^{\ast}))$, we get
$$
\int_{\RN} \left(F(u_{n,\e}^R-u_\e^R)-\frac{1}{2^{\ast}}|u_{n,\e}^R-u_\e^R|^{2^\ast}\right)=o(1).
$$
Thus, using the Sobolev's inequality, we get that for any $n$ large,
{\allowdisplaybreaks
\begin{align*}
D_\e & \ge  \G_\e(u_{n,\e}^R) \\
&\ge \G_\e(u_\e^R) + \frac 12\| w_{n,\e}^R-w_\e^R\|^2 - \frac{1}{2^*}\|w_{n,\e}^R-w_\e^R\|_{2^*}^{2^*} + o(1)\\
&\ge\G_\e(u_\e^R)+\frac 12\|w_{n,\e}^R-w_\e^R \|^{2}-\frac{1}{2^*}S^{-N/(N-2)}\|w_{n,\e}^R-w_\e^R\|^{2^*}+ o(1)\\
&=\G_\e(u_\e^R)+\|w_{n,\e}^R-w_\e^R\|^{2}\Big(\frac 12-\frac{1}{2^*}S^{-N/(N-2)}\|w_{n,\e}^R-w_\e^R\|^{4/(N-2)}\Big )+o(1)\\
&\ge\G_\e(u_\e^R)+o(1).
\end{align*}
}%
Letting $n\to+\iy$, we get that $\G_\e(u_\e^R)\le D_\e$, that is, $u_\e^R\in \G_\e^{D_\e}$.

\vskip0.1in

\noindent {\bf Step 2.} We claim that for $d\in(0,\frac{1}{2}S^{\frac{N}{4}})$, $\{u_\e^R\}$ is bounded in $L^\infty(\RN)$ uniformly for large $R>0$ .
To the contrary, we assume that there exist $R_0,R_n>0$ satisfying $R_n>R_0$ and $\lim_{n\rg \iy}\|u_\e^{R_n}\|_{L^{\iy}(\RN)}=\iy$, where $\{u_\e^{R_n}\}\subset X_\e^d\cap\G_\e^{D_{\e}}$.
By the definition of $X_\e^d$, there exist $\{U_n\}\subset S_m$ and $\{y_n\}\subset\mathcal{M}^\bb$ with
$$
\|u_\e^{R_n}-\varphi_\e(\cdot-\frac{y_n}{\e})U_n(\cdot-\frac{y_n}{\e})\|_\e\le d.
$$
Since $S_m$ is compact, up to a subsequence, there exists $U\in S_m$ such that $U_n\rg U$ in $H^1(\RN)$. Then
\begin{eqnarray*}
& &\|u_\e^{R_n}-\varphi_\e(\cdot-\frac{y_n}{\e})U(\cdot-\frac{y_n}{\e})\|_\e\\
& & \le  d+\|\varphi_\e(\cdot-\frac{y_n}{\e})U_n(\cdot-\frac{y_n}{\e})-\varphi_\e(\cdot-\frac{y_n}{\e})U(\cdot-\frac{y_n}{\e})\|_\e\\
&  & =:d+I_n.
\end{eqnarray*}
It follows from $U_n\rg U$ in $H^1(\RN)$ that $\lim_{n\rg\iy}I_n=0$, which implies that
$
\limsup_{n\rg\infty}\|u_\e^{R_n}-\varphi_\e(\cdot-\frac{y_n}{\e})U(\cdot-\frac{y_n}{\e})\|_\e\le \frac{3d}{2}.
$
Then we write
$$
u_\e^{R_n}=v_\e^{R_n}+w_\e^{R_n}\ \ \mbox{and}\ \ \|v_\e^{R_n}\|_\e\le 2d,
$$
for $n$ large enough, where $v_\e^{R_n}\in H_\e$ and $w_\e^{R_n}=\varphi_\e(\cdot-\frac{y_n}{\e})U(\cdot-\frac{y_n}{\e})$. Meanwhile, since $f(t)=0$ for $t\le 0$, it follows form the Maximum Principle that $u_\e^{R_n}>0$ in $B(0,\frac {R_n}{\e})$ and $-\DD u_\e^{R_n}+u_\e^{R_n}\le f(u_\e^{R_n})$ in $B(0,\frac {R_n}{\e})$. Due to elliptic estimates, $u_\e^{R_n}\in C^{1,\alpha}(B(0,\frac {R_n}{\e}))$ for some $\alpha\in (0,1)$. We extend $u_\e^{R_n}\in H_0^1(B(0,\frac {R_n}{\e}))$ to $u_\e^{R_n}\in H^1(\RN)$ by zero outside $B(0,\frac {R_n}{\e})$. Then, there exists $z_\e^n\in B(0,\frac {R_n}{\e})$ such that $u_\e^{R_n}(z_\e^n)=\max_{x\in \RN}u_\e^{R_n}(x)$ for $n = 1,2,\cdots.$ We define
\[w_n(x) \equiv \frac{1}{l_n}u_\e^{R_n}(l_n^{-\frac{2}{N-2}}x+z_\e^n), \ \ w^0_n(x) \equiv \frac{1}{l_n}w_\e^{R_n}(l_n^{-\frac{2}{N-2}}x+z_\e^n) \ \ \textrm{with} \ \ l_n \equiv u_\e^{R_n}(z_\e^n).\]
Then, $w_n$ satisfies
{\allowdisplaybreaks
\begin{align*}
&-\DD w_n +(l_n)^{-\frac{4}{N-2}}w_n -(l_n)^{-\frac{N+2}{N-2}}f(l_nw_n)\\
&=-4(l_n)^{-\frac{4}{N-2}}\left(\int_{\RN}\chi_\e(x) (u_\e^{R_n})^2dx-1\right)_+\chi_\e(l_n^{-\frac{2}{N-2}}x+z_\e^n)w_n\ \ \ \textrm { in  } \RN.
\end{align*}
}%
Note that $\{\int_{\RN}\chi_\e (u_\e^{R_n})^2\}$ and $\{\|u_\e^{R_n}\|_\e\}$ are uniformly bounded for $n$. Then, we deduce from (F2) and elliptic estimates that $w_n$ converges locally uniformly to
the unique radial positive solution $w_0$ of
\[ \DD w_0 + w_0^{\frac{N+2}{N-2}} = 0,   \ \ w_0(0)=\max_{x\in\RN}w_0(x) = 1 \ \ \textrm { in  } \RN.\]
It is  well known that $w_0(x)= \Big (\frac{N(N-2)}{N(N-2)+|x|^2}\Big )^{\frac{N-2}{2}}$; thus $w_0 \in L^{2^*}(\RN)$ and $|\nabla w_0| \in L^2(\RN).$ Note that $\|w_n-w_n^0\|_{2^*} = \|u_\e^{R_n}-w_\e^{R_n}\|_{2^*}=\|v_\e^{R_n}\|_{2^*}\le \frac{2d}{\sqrt{S}}$ for $n$ large enough. Then for any fixed $R>0$, we have
$\|w_n-w_n^0\|_{L^{2^\ast}(B(0,R))}\le \frac{2d}{\sqrt{S}}$. Recall that $\|w_n^0\|_{L^\iy}\le\|U\|_{L^\iy}/l_n\to 0$ as $n\to\iy$, letting $n\to\iy$ we deduce that
$\|w_0\|_{L^{2^\ast}(B(0,R))}\le \frac{2d}{\sqrt{S}}$. This means that $\|w_0\|_{2^*}\le \frac{2d}{\sqrt{S}};$ but this is impossible since $d\in(0,\frac{1}{2}S^{\frac{N}{4}})$. Thus $\{u_\e^R\}$ is bounded in $L^\infty(\RN)$ uniformly for $R>R_0$.

\vskip0.1in

\noindent {\bf Step 3.} We claim that $u_\e^R\rg u_\e$ strongly in $H_\e$, where $\G_\e'(u_\e)=0, u_\e\in X_\e^d\cap\G_\e^{D_\e}$. By Step 2 and elliptic estimates (see \cite{Tru}), we see that there exists $C>0$ (independent $R$) such that for any $B(y,2)\subset B(0,\frac{R}{\e})$, $\sup_{B(y,1)}u_\e^R\le C\|u_\e^R\|_{L^2(B(y,2))}$. Thus, by $Q_\e(u_\e^R)$ is uniformly bounded for $R>R_0$, we see that there exists $C>0$ (independent $R$) such that
\be\lab{ss1}
0<u_\e^R\le C\e^\frac{\mu}{2}\ \mbox{for}\ \ |y|\ge 2+\frac{R_0}{\e}, R>R_0.
\ee
Since $u_\e^R\in X_\e^d$, we can assume $u_\e^R\rightharpoonup u_\e$ weakly in $H_\e^1(\RN)$ as $R\rg\iy$. Next, we shall prove that $u_\e^R\rg u_\e$ strongly in $H_\e^1(\RN)$ as $R\rg\iy$. First, we claim that the sequence $\{u_\e^R\}_{R>R_0}$ has exponential decay at infinity. By (\ref{ss1}) and $(F1)$, for sufficiently small and fixed $\e>0$, we have $|f(u_\e^R)|\le \frac 12 u_\e^R$ for $|y|\ge 2+\frac{R_0}{\e}, R>R_0$. It follows from the Maximum Principle that $0<u_\e^R(y)\le C\exp(-\frac{1}{2}|y|)$ for $|y|\ge 2+\frac{R_0}{\e}, R>R_0$. Therefore,  by Step 2 there exists $C>0$ (independent $R$) such that
\be\lab{ss2}
0<u_\e^R(y)\le C\exp(-\frac{1}{2}|y|)\ \mbox{for}\ y\in\RN, R>R_0.
\ee
Second, we claim that
\begin{eqnarray}
\lab{ss3}
\lim_{\dd\rg \iy}\int_{\RN\setminus B(0,\dd)}|\na u_\e^R|^2+V_\e|u_\e^R|^2=0,
\end{eqnarray}
uniformly for $R>R_0$. Choosing a cutoff function $\phi_\dd\in C^\iy(\RN), 0\le\phi_\dd\le1, |\na \phi_\dd|\le\frac{2}{\dd}$ and $\phi_\dd(y)=0, |y|\le \dd, \phi_\dd(y)=1, |y|\ge 2\dd$, it follows from  $\langle\G_\e'(u_\e^R),\phi_\dd u_\e^R\rangle=0$ that
\begin{eqnarray}
\lab{ss4}
&\ &\int_{\RN\setminus B(0,2\dd)}|\na u_\e^R|^2+V_\e|u_\e^R|^2\nonumber\\
&\ &\le\frac{1}{\dd}\int_{\RN\setminus B(0,\dd)}|\na u_\e^R|^2+|u_\e^R|^2+\int_{\RN\setminus B(0,\dd)}|f(u_\e^R)u_\e^R|.
\end{eqnarray}
Thus, (\ref{ss3}) immediately follows form (\ref{ss2}),(\ref{ss4}) and the fact that $\{\|u_\e^R\|_\e\}_{R>R_0}$ is uniformly bounded. Third, we shall prove that $\G_\e'(u_\e)=0$ in $H_\e(\RN)$ and $u_\e^R\rg u_\e$ strongly in $H_\e(\RN)$. By (\ref{ss3}) we see that $\lim_{\dd\rg\iy}\|u_\e^R\|_{\RN\setminus B(0,\dd)}=0$. Then by $u_\e^R\rg u_\e$ weakly in $H_\e(\RN)$, we get that $u_\e^R\rg u_\e$ strongly in $L^q(\RN)(2\le q<2^{\ast})$. Since $\{u_\e^R\}$ is bounded in $L^\infty(\RN)$ uniformly for $R>R_0$, we see that $\|u_\e\|_\iy<\iy$. So $u_\e^R\rg u_\e$ strongly in $L^q(\RN)(2\le q\le2^{\ast})$. Therefore, by a standard way we can prove the claim.

Finally, since $S_m$ is compact, it is easy to see that $0\not\in X_\e^d$ for $d>0$ small enough. Thus, $u_\e\not\equiv0$. This completes the proof.
\ep

Now, we prove Theorem \ref{Theorem 1}. We  start with the following Lemma \ref{B-K} due to Br$\acute{e}$zis and Kato \cite{B-K} and Lemma \ref{Tru} due to Gilbarg and Trudinger \cite{Tru}.

\bl\lab{B-K}\cite{B-K}
Let $\Omega\subset\RN$ and $q\in L^{\frac{N}{2}}(\RN)$ be a nonnegative function. Then for every $\e>0$, there exists a constant $\si (\e,q)>0$ such that
$$
\int_{\Omega}q(x)u^2\le \e\int_{\Omega} |\na u|^2+\si (\e,q)\int_{\Omega}u^2,\ \ \mbox{for all}\ \ u\in H^1 (\Omega).
$$
\el

\bl\lab{Tru}\cite{Tru}
Suppose that $t>N,g\in L^{\frac{t}{2}}(\om)$ and $u\in H_0^1(\om)$ satisfies (in the weak sense)
$$
-\DD u+u\le g(x),
$$
where $\om$ is an open subset of $\RN$. Then for any ball $B_{2R}(y)\subset \om$,
$$
\sup_{B_R(y)}u\le C\left(\|u^+\|_{L^2(B_{2R}(y))}+\|g\|_{L^{\frac{t}{2}}(B_{2R}(y))}\right),
$$
where $C$ only depends on $N,t$ and $R$.
\el

\noindent{\bf Completion of the proof for Theorem \ref{Theorem 1}}
\bp
By Proposition \ref{prop2.6}, there exist $d>0$ and $\e_0>0$, such that $\G_\e$ has a nontrivial critical point $u_\e\in X_\e^d\cap\G_\e^{D_\e}$ for $\e\in(0,\e_0)$.
\vskip0.1in

\noindent {\bf Step 1.} We claim that for $d>0$ small, there exists $\rho>0$ (independent of $\e$) such that $\|u_\e\|_\iy\ge\rho$ for $\e\in(0,\e_0)$ and $u_\e>0$ in $\RN$. Since $f(t)=0$ for $t\le 0$, we see that $u_\e\ge 0$  and $-\DD u_\e+u_\e\le f(u_\e)$ in $\RN$. In the following, we use the Moser iteration technique (see \cite{Tru}) to prove that there exists $C>0$ such that
\be\lab{liy}
\|u_\e\|_\iy<C,\ \mbox{uniformly for}\ \e\in(0,\e_0).
\ee
If (\ref{liy}) is ture, then it follows from weak Harnark inequality (see \cite{Tru}) that $u_\e>0$ in $\RN$. Thus, from $(V1)$ and $(F1)$, it is easy to see that $\inf_{\e\in(0,\e_0)}\|u_\e\|_\iy>0$.

To the contrary, we assume that (\ref{liy}) is false, i.e., there exist  $\e_n,d_n>0$ satisfying $\lim\limits_{n\rg \iy}d_n=\lim\limits_{n\rg \iy}\e_n$ $=$ $0$ and $\lim\limits_{n\rg \iy}\|u_{\e_n}\|_{L^{\iy}(\RN)}=\iy$, where $\{u_{\e_n}\}\subset X_{\e_n}^{d_n}\cap\G_{\e_n}^{D_{\e_n}}$. By the definition of $X_{\e_n}^{d_n}$, there exist $\{U_n\}\subset S_m$ and $\{y_n\}\subset \mathcal{M}_{\e_n}^\bb$ with
$$
\|u_{\e_n}-(\vp_{\e_n}U_n)(\cdot-y_n)\|_{\e_n}\le d_n.
$$
Since $S_m$ is compact, up to a subsequence, there exists $U\in S_m$ such that $U_n\rg U$ in $H^1(\RN)$. Then
\begin{eqnarray*}
& &\|u_{\e_n}-(\vp_{\e_n}U)(\cdot-y_n)\|_{\e_n}\\
& & \le  d_n+\|(\vp_{\e_n}U_n)(\cdot-y_n)-(\vp_{\e_n}U)(\cdot-y_n)\|_{\e_n}\\
&  & =: d_n+I_n.
\end{eqnarray*}
It follows from $U_n\rg U$ in $H^1(\RN)$ that $I_n\rg 0$ as $n\rg\iy$,
which implies that
$
\lim_{n\rg\infty}\|u_{\e_n}-(\vp_{\e_n}U)(\cdot-y_n)\|_{\e_n}=0.
$
By Sobolev's embedding theorem, for any $\mu>0$, there exists $n_0\in \mathbb{N}$ such that
\be\lab{qq10}
\|u_{\e_n}^{2^{\ast}-2}-(\vp_{\e_n}U)^{2^{\ast}-2}(\cdot-y_n)\|_{L^{\frac{N}{2}}(\RN)}\le \mu\ \ \mbox{for}\ \ n\ge n_0.
\ee
For convenience, we write $n$ for $\e_n$. For any $k\in \mathbb{N}$ and $p>0$, consider $A_k=\{x\in\RN:u_n\le k\}, B_k=\RN\setminus A_k$ and define $v_k$ by
$$
v_k=u_n^{2p+1}\ \mbox{in}\ \ A_k,\ \ v_k=k^{2p}u_n,\ \mbox{in}\ \ B_k.
$$
Thus, $v_k\in H^1(\RN), 0\le v_k\le u_n^{2p+1}$ and
$$
\na v_k=(2p+1)u_n^{2p}\na u_n\ \mbox{in}\ \ A_k,\ \ \na v_k=k^{2p}\na u_n,\ \mbox{in}\ \ B_k.
$$
So, using $v_k$ as a test function, we have $\int_{\RN}\na u_n\na v_k\le\int_{\RN}f(u_n)v_k$, i.e.
$$
(2p+1)\int_{A_k}u_n^{2p}|\na u_n|^2+k^{2p}\int_{B_k}|\na u_n|^2=\int_{\RN}f(u_n)v_k.
$$
By $(F1)$-$(F2)$, there exists $C>0$ such that $f(t)\le Ct+2t^{2^{\ast}-1}$ for all $t>0$. Thus, we get
\be\lab{q7}
(2p+1)\int_{A_k}u_n^{2p}|\na u_n|^2+k^{2p}\int_{B_k}|\na u_n|^2\le\int_{\RN}(C+2u_n^{2^{\ast}-2})u_nv_k.
\ee
We define
$$
w_k \equiv u_n^{p+1}\ \mbox{in}\ \ A_k,\ \ w_k \equiv k^pu_n\ \mbox{in}\ \ B_k.
$$
Then $w_k^2=u_nv_k\le u_n^{2p+2}$ and
$$
\na w_k=(p+1)u_n^p\na u_n\ \mbox{in}\ \ A_k,\ \ \na w_k=k^p\na u_n,\ \mbox{in}\ \ B_k.
$$
Thus we get
\be\lab{q8}
\int_{\RN}|\na w_k|^2=(p+1)^2\int_{A_k}u_n^{2p}|\na u_n|^2+k^{2p}\int_{B_k}|\na u_n|^2.
\ee
Combining (\ref{q7}) and (\ref{q8}), we see that
\be\lab{q9}
\frac{2p+1}{(p+1)^2}\int_{\RN}|\na w_k|^2\le\int_{\RN}(C+2u_n^{2^{\ast}-2})w_k^2.
\ee
For any $\mu>0$ given above, it follows from Lemma \ref{B-K} and (\ref{qq10}) that
\begin{align*}
\lab{q11}
\int_{\RN}u_n^{2^{\ast}-2}w_k^2&=\int_{\RN}u_0^{2^{\ast}-2}w_k^2+\int_{\RN}\left(u_n^{2^{\ast}-2}-u_0^{2^{\ast}-2}\right)w_k^2\\
&\le\mu\int_{\RN} |\na w_k|^2+\si (\mu,U)\int_{\RN}w_k^2+\mu\|w_k\|_{L^{2^{\ast}}(\RN)}^2\\
&\le\frac{1+S}{S}\mu\int_{\RN} |\na w_k|^2+\si (\mu,U)\int_{\RN}w_k^2,
\end{align*}
where $S$ is the best Sobolev's constant. Then, choosing $\mu=\frac{S(2p+1)}{4(S+1)(p+1)^2}$, we get from (\ref{q9}) that
\be\lab{q12}
\int_{\RN}|\na w_k|^2\le \ti{C_p}\int_{\RN}|w_k|^2,
\ee
where $\ti{C_p}=\frac{2(p+1)^2}{2p+1}\left(C+2\si(\mu,U)\right)$. That is, for any $p>0$, there exists $n_p\in \mathbb{N}$ such that
\be\lab{q13}
\int_{\RN}|\na w_k|^2\le \ti{C_p}\int_{\RN}|w_k|^2 \ \mbox{for any}\ k>0\ \mbox{and}\ n\ge n_p.
\ee
If $u_n\in L^{2(p+1)}(\RN)$ for some $p\ge 2$, by Sobolev's embedding theorem and (\ref{q13}), we have
$$\left(\int_{A_k}w_k^{2^{\ast}}\right)^{\frac{2}{2^{\ast}}}\le S\ti{C_p}\int_{\RN}w_k^2.$$
Therefore, for any $p\ge 2$ satisfying $\{u_n\}_{n=1}^\iy\subset L^{2(p+1)}(\RN)$, there exists $\ti{C_p}>0,n_p\in \mathbb{N}$ such that
$$\left(\int_{A_k}u_n^{2^{\ast}(p+1)}\right)^{\frac{2}{2^{\ast}}}\le S\ti{C_p}\int_{\RN}u_n^{2(p+1)},$$
for any $k>0$ and $n\ge n_p$. Now, let $k\rg\iy$, we have
\be\lab{q14}
\left(\int_{\RN}u_n^{2^{\ast}(p+1)}\right)^{\frac{2}{2^{\ast}}}\le S\ti{C_p}\int_{\RN}u_n^{2(p+1)}.
\ee
That is, for any $p\ge 2$ satisfying $u_n\subset L^{2(p+1)}(\RN)$ for $n$ large enough, there exist $C_p>0,n_p\in \mathbb{N}$, such that for $n\ge n_p$, $u_n\in L^{2^{\ast}(p+1)}(\RN)$ and
\be\lab{q15}
\|u_n\|_{L^{2^{\ast}(p+1)}(\RN)}\le C_p\|u_n\|_{L^{2(p+1)}(\RN)}.
\ee
In the following, we will use an iteration argument. Let $p_1$ be a positive constant such that $2(p_1+1)=2^{\ast}$. Noting that $\{u_n\}\subset L^{2^{\ast}}(\RN)$, by (\ref{q15}), there exist $C_1>0,n_1\in \mathbb{N}$, such that for $n\ge n_1$, $u_n\in L^{2^{\ast}(p_1+1)}(\RN)$ and
$$
\|u_n\|_{L^{2^{\ast}(p_1+1)}(\RN)}\le C_1\|u_n\|_{L^{2(p_1+1)}(\RN)}.
$$
Choosing $p_2$ satisfying  $2(p_2+1)=2^{\ast}(p_1+1)$, we see that $p_2>p_1$ and $u_n\subset L^{2(p_2+1)}(\RN)$ for $n\ge n_1$. Thus, by (\ref{q15}), there exist $C_2>0,n_2\in \mathbb{N}$ and $n_2\ge n_1$, such that for $n\ge n_2$, $u_n\in L^{2^{\ast}(p_2+1)}(\RN)$ and
$$
\|u_n\|_{L^{2^{\ast}(p_2+1)}(\RN)}\le C_2\|u_n\|_{L^{2(p_2+1)}(\RN)}.
$$
Continuing with this iteration, we get a consequence $\{C_k\}$ and two increasing sequences $\{n_k\}$ and $\{p_k\}$, where $2(p_{k+1}+1)=2^{\ast}(p_k+1)$, such that for $n\ge n_k$, $u_n\in L^{2^{\ast}(p_k+1)}(\RN)$ and
$$
\|u_n\|_{L^{2^{\ast}(p_k+1)}(\RN)}\le C_k\|u_n\|_{L^{2(p_k+1)}(\RN)}.
$$
Obviously, $p_k=\left(\frac{N}{N-2}\right)^{k-1}2^{\ast}-1$. It follows that, for any $p\ge2$ there exist $C_p>0$ and $n_p\in \mathbb{N}$ such that for $n\ge n_p$,
\be\lab{q16}
\|u_n\|_{L^p(\RN)}\le C_p\|u_n\|_{L^{2^{\ast}}(\RN)}.
\ee
On the other hand, by $(F1)$-$(F2)$, $|f(t)|\le C(|t|+|t|^{2^{\ast}-1})$ for all $t$ and some $C>0$. By (\ref{q16}), there exists $n_0\in \mathbb{N}$, such that for $n\ge n_0$,
\be\lab{q17}
\|f(u_n)\|_{L^N(\RN)}\le C_0\left(\|u_n\|_{L^{2^{\ast}}(\RN)}+\|u_n\|_{L^{2^{\ast}}(\RN)}^{2^{\ast}-1}\right)
\ee
for some $C_0>0$. By Lemma \ref{Tru}, there exist $n_N\in\mathbb{N}$ such that, for any $y\in \RN$,
\be\lab{q18}
\sup_{B_1(y)}u_n\le C\left(\|u_n\|_{L^2(B_{2}(y))}+\|f(u_n)\|_{L^N(B_{2}(y))}\right),\ \mbox{for}\ n\ge n_N,
\ee
where $C$ only depends on $N$. Obviously, $\{\|u_n\|_2\}, \{\|u_n\|_{2^\ast}\}$ are bounded uniformly for $n$. Then we see that $\sup_{n\ge n_N}\|u_n\|_{L^\iy(\RN)}<\iy$, which is a contradiction. Therefore, the claim (\ref{liy}) is concluded.

\vskip0.1in

\noindent {\bf Step 2.} There exist $M>0$ (independent of $\e$) and $y_\e\in\RN$, such that
\be\lab{ss5}
0<w_\e(y)\le M\exp\left(-\frac{|y|}{2}\right),\ \mbox{for}\ y\in\RN, \e\in(0,\e_0),
\ee
where $w_\e(y)=u_\e(y+y_\e)$. By Proposition \ref{prop2.4}, for small $d>0$ there exist $\{y_\e\}\subset\RN, x\in \mathcal{M}, U\in S_m$ such that
$$
\lim_{\e\rg 0}|\e y_\e-x|=0\ \mbox{and}\ \lim_{\e\rg 0}\|u_\e-U(\cdot-y_\e)\|_\e=0.
$$
Then for any $\sigma>0$, there exists $R>0$ (independent of $\e$), such that
$$
\sup_{\e\in(0,\e_0)}\int_{\RN\setminus B(0,R)}w_\e^2\le\sigma.
$$
Moreover, since $\G_\e'(u_\e)=0$ and $\{u_\e\}$ is bounded in $L^\infty(\RN)$ uniformly for $\e\in(0,\e_0)$, there exists $C>0$ (independent of $\e$), such that $-\DD w_\e\le Cw_\e$ in $\RN$. Then by elliptic estimates (see \cite{Tru}), there exists $C>0$ (independent of $\e$), such that $\sup_{B(y,1)}w_\e\le C\|w_\e\|_{L^2(B(y,2))}$ for any $y\in \RN$. Then, $0<w_\e(y)\le C\sigma^{\frac12}$ for $ \e\in(0,\e_0),|y|\ge R+2$. Thus, the claim can be proved by Maximum Principle.
\vskip0.1in

\noindent {\bf Step 3.} We claim that $Q_\e(u_\e)=0$ for small $\e>0$. Since $\lim_{\e\rg 0}\e y_\e=x\in \mathcal{M}$,  $\e y_\e\in \mathcal{M}^{5\delta}$ for small $\e$ and there is $C>0$ such that $|y|\le Cdist(y,\mathcal{M}^{5\delta})$ for $y\in \RN\setminus O$. Thus there is $C>0$ such that
{\allowdisplaybreaks
\begin{align*}
\int_{\RN}\chi_\e(u_\e)^2&\le C\e^{-\mu}\int_{\RN\setminus O_\e}\exp{(-\frac{1}{2}|x-y_\e|)}dx\\
&=C\e^{-\mu-N}\int_{\RN\setminus O}\exp{(-\frac{1}{2\e}|y-\e y_\e|)}dy\\
&\le C\e^{-\mu-N}\int_{\RN\setminus O}\exp{(-\frac{1}{2\e}dist(y,\mathcal{M}^{5\delta}))}dy\\
&\le C\e^{-\mu-N}\int_{\RN\setminus B(0,\delta)}\exp{(-\frac{1}{C\e}|y|)}dy\\
&\longrightarrow0,\ \ \mbox{as}\ \ \e\rg 0,
\end{align*}
}%
i.e., $Q_\e(u_\e)=0$ for small $\e>0$. Therefore, $u_\e$ is a critical point of $P_\e$ and a solution of (\ref{q2}).
\vskip0.1in

\noindent {\bf Step 4.} We shall prove that there exists $z_0,x_\e\in \RN$, such that $\max_{\RN}v_\e=v_\e(x_\e)$, $\lim_{\e\rg 0}dist(x_\e,\mathcal{M})=0$ and $\lim_{\e\rg 0}\|v_\e(\e\cdot+x_\e)-U(\cdot+z_0)\|_\e=0$. Assume that $z_\e\in \RN$ such that $\|w_\e\|_\iy=w_\e(z_\e)$, then by Step 1 and (\ref{ss5}), we see that $\{z_\e\}\subset\RN$ is bounded. Up to a subsequence, we can assume that $z_\e\rg z_0$ as $\e\rg 0$. Let $\ti{x_\e}=y_\e+z_\e$, then $\max_{\RN}u_\e=u_\e(\ti{x_\e})$. Thus, let $x_\e=\e y_\e+\e z_\e$, we see that $\max_{\RN}v_\e=v_\e(x_\e)$ and $\lim_{\e\rg 0}x_\e=\lim_{\e\rg 0}\e y_\e=x\in \mathcal{M}$. Finally, it is easy to check that $v_\e(\e\cdot+x_\e)\rg U(\cdot+z_0)$ strongly in $H_\e(\RN)$ as $\e\rg 0$. This completes the proof.
\ep

\end{document}